# GENERALIZED HERMITE-HADAMARD TYPE INEQUALITIES INVOLVING LOCAL FRACTIONAL INTEGRALS

HUIXIA MO

School of Science, Beijing University of Posts and Telecommunications, Beijing, 100876, P. R. China
Email: huixmo@bupt.edu.cn

In the paper, two new identities involving the local fractional integrals have been established. Using these two identities, we obtain some generalized Hermite-Hadamard type integral inequalities for the local differentiable generalized convex functions.

## 1. INTRODUCTION

The convex function is an important concept in the class mathematical analysis course. And many important inequalities are established for the class of convex functions. For example, the Hermite-Hadamard's inequality is one of the best known results in the literature, which can be stated as follows.

**Hermite-Hadamard's inequality**[1]: Let $f$ be a convex function on $[a,b]$ with $a < b$. If $f$ is integral on $[a,b]$, then

$$f\left(\frac{a+b}{2}\right) \leq \frac{1}{b-a}\int_a^b f(x)dx \leq \frac{f(a)+f(b)}{2}. \tag{1.1}$$

Hadamard's inequality for convex functions has received renewed attention in recent years and a remarkable variety of refinements and generalizations have been found(see [2-6] and so on). In [2], Dragomir and Agarwal proved the following results which connected with the right part of (1.1).

**THEOREM A.** [2] Let $I \subseteq R$ be an interval. Suppose that $f : I^0 \subseteq R \to R$ ($I^0$ is the interior of $I$) is differentiable. If $f' \in L[a,b]$ and $|f'|$ is convex on $[a,b]$ for $a,b \in I^0$ with $a < b$, then

$$\left|\frac{f(a)+f(b)}{2} - \frac{1}{b-a}\int_a^b f(x)dx\right| \leq \frac{b-a}{8}(|f'(a)| + |f'(b)|). \tag{1.2}$$

**THEOREM B.** [2] Let $f : I^0 \subseteq R \to R$ be a differentiable mapping on $I^0$. If $f' \in L[a,b]$ and $|f'|^{p/(p-1)}$ is convex on $[a,b]$ for $p > 1$ and $a,b \in I^0$ with $a < b$, then

$$\left|\frac{f(a)+f(b)}{2} - \frac{1}{b-a}\int_a^b f(x)dx\right|$$
$$\leq \frac{b-a}{2(1+p)^{1/p}}\left[\frac{|f'(a)|^{p/(p-1)} + |f'(b)|^{p/(p-1)}}{2}\right]^{(p-1)/p}. \tag{1.3}$$

In [5], the above inequalities were generalized as follows.

**THEOREM C.** [5] Let $f : I^0 \subseteq R \to R$ be a differentiable mapping on $I^0$. If $f' \in L[a,b]$ and $|f'|^q$ is convex on $[a,b]$ for $q \geq 1$ and $a,b \in I^0$ with $a < b$, then the following equality holds:

$$\left| \frac{f(a)+f(b)}{2} - \frac{1}{b-a}\int_a^b f(x)dx \right| \leq \frac{b-a}{4}\left[ \frac{|f'(a)|^q + |f'(b)|^q}{2} \right]^{1/q},$$

$$\left| f\left(\frac{a+b}{2}\right) - \frac{1}{b-a}\int_a^b f(x)dx \right| \leq \frac{b-a}{4}\left[ \frac{|f'(a)|^q + |f'(b)|^q}{2} \right]^{1/q}.$$

In recent years, the fractal theory has received significantly remarkable attention [7-11]. In [9], the authors introduce the generalized convex function and establish the generalized Hermite-Hadamard's inequality on fractal space.

**THEOREM D. [9] (Generalized Hermite-Hadmard's inequality)** Let $f(x) \in {_aI_b^{(\alpha)}}$ be a generalized convex function on $[a,b]$ with $a < b$. Then

$$f\left(\frac{a+b}{2}\right) \leq \frac{\Gamma(1+\alpha)}{(b-a)^\alpha} {_aI_b^{(\alpha)}} f(x) \leq \frac{f(a)+f(b)}{2^\alpha}.$$

Inspired by these investigations, in the paper we will establish the generalized Hermite-Hadamard type integral inequalities for the local differentiable convex functions on fractal space.

## 2. PRELIMINARIES

Recall the set $R^\alpha$ of real line numbers and use the Gao-Yang-Kang's idea to describe the definitions of the local fractional derivative and local fractional integral, see [10,11] and so on.

Recently, the theory of Yang's fractional sets [10] was introduced as follows.

For $0 < \alpha \leq 1$, we have the following $\alpha$-type set of element sets:

$Z^\alpha$: The $\alpha$-type set of the integer is defined as the set $\{0^\alpha, \pm 1^\alpha, \pm 2^\alpha, \cdots, \pm n^\alpha, \cdots\}$.

$Q^\alpha$: The $\alpha$-type set of the rational numbers is defined as the set $\{m^\alpha = (p/q)^\alpha: p, q \in Z, q \neq 0\}$.

$J^\alpha$: The $\alpha$-type set of the irrational numbers is defined as the set $\{m^\alpha \neq (p/q)^\alpha: p, q \in Z, q \neq 0\}$.

$R^\alpha$: The $\alpha$-type set of the real line numbers is defined as the set $R^\alpha = Q^\alpha \cup J^\alpha$.

If $a^\alpha, b^\alpha$ and $c^\alpha$ belong to the set $R^\alpha$ of real line numbers, then

(1) $a^\alpha + b^\alpha$ and $a^\alpha b^\alpha$ belong to the set $R^\alpha$;                (2) $a^\alpha + b^\alpha = b^\alpha + a^\alpha = (a+b)^\alpha = (b+a)^\alpha$;

(3) $a^\alpha + (b^\alpha + c^\alpha) = (a+b)^\alpha + c^\alpha$;                (4) $a^\alpha b^\alpha = b^\alpha a^\alpha = (ab)^\alpha = (ba)^\alpha$;

(5) $a^\alpha(b^\alpha c^\alpha) = (a^\alpha b^\alpha)c^\alpha$;                (6) $a^\alpha(b^\alpha + c^\alpha) = a^\alpha b^\alpha + a^\alpha c^\alpha$;

(7) $a^\alpha + 0^\alpha = 0^\alpha + a^\alpha = a^\alpha$ and $a^\alpha 1^\alpha = 1^\alpha a^\alpha = a^\alpha$.

The definitions of the local fractional derivative and local fractional can be given as follows.

**DEFINITION 2.1.**[10] A non-differentiable function $f: R \to R^\alpha, x \to f(x)$ is called to be local fractional continuous at $x_0$, if for any $\varepsilon > 0$, there exists $\delta > 0$, such that

$$|f(x) - f(x_0)| < \varepsilon^\alpha$$

holds for $|x-x_0|<\delta$, where $\varepsilon,\delta \in R$. If $f(x)$ is local fractional continuous on the interval $(a,b)$, we denote $f(x) \in C_\alpha(a,b)$.

**DEFINITION 2.2.**[10] The local fractional derivative of $f(x)$ of order $\alpha$ at $x = x_0$ is defined by

$$f^{(\alpha)}(x_0) = \frac{d^\alpha f(x)}{dx^\alpha}\bigg|_{x=x_0} = \lim_{x \to x_0} \frac{\Delta^\alpha(f(x)-f(x_0))}{(x-x_0)^\alpha},$$

where $\Delta^\alpha(f(x)-f(x_0)) \cong \Gamma(1+\alpha)(f(x)-f(x_0))$.

If there exists $f^{((k+1)\alpha)}(x) = \overbrace{D_x^\alpha \cdots D_x^\alpha}^{k+1 \text{ times}} f(x)$ for any $x \in I \subseteq R$, then we denoted $f \in D_{(k+1)\alpha}(I)$, where $k = 0,1,2,\ldots$.

**DEFINITION 2.3.**[10] Let $f \in C_\alpha[a,b]$. Then the local fractional integral is defined by,

$$_aI_b^{(\alpha)} f(x) = \frac{1}{\Gamma(1+a)} \int_a^b f(t)(dt)^\alpha = \frac{1}{\Gamma(1+a)} \lim_{\Delta t \to 0} \sum_{j=0}^{N-1} f(t_j)(\Delta t_j)^\alpha,$$

with $\Delta t_j = t_{j+1}-t_j$ and $\Delta t = \max\{\Delta t_1, \Delta t_2, \Delta t_j,\ldots\}$, where $[t_j, t_j+1]$, $j = 0,\ldots,N-1$ and $t_0 = a < t_1 < \cdots < t_i < \cdots < t_{N-1} < t_N = b$ is a partition of the interval $[a,b]$.

Here, it follows that $_aI_a^{(\alpha)} f(x) = 0$ if $a = b$ and $_aI_b^{(\alpha)} f(x) = -_bI_a^{(\alpha)} f(x)$ if $a < b$. If for any $x \in [a,b]$, there exists $_aI_x^{(\alpha)} f(x)$, then we denoted by $f(x) \in I_x^{(\alpha)}[a,b]$.

**DEFINITION 2.4.** [9] Let $f : I \subseteq R \to R^\alpha$. For any $x_1, x_2 \in I$ and $\lambda \in [0,1]$, if the following inequality

$$f(\lambda x_1 + (1-\lambda)x_2) \leq \lambda^\alpha f(x_1) + (1-\lambda)^\alpha f(x_2)$$

holds, then $f$ is called a generalized convex function on $I$.

Here are two basic examples of generalized convex functions:

(1) $f(x) = x^{\alpha p}$, $x \geq 0$, $p > 1$;

(2) $f(x) = E_\alpha(x^\alpha), x \in R$, where $E_\alpha(x^\alpha) = \sum_{k=0}^{\infty} \frac{x^{\alpha k}}{\Gamma(1+k\alpha)}$ is the Mittag-Leffer function.

**LEMMA 2.1.**[10]

(1) **(Local fractional integration is anti-differentiation)** Suppose that $f(x) = g^{(\alpha)}(x) \in C_\alpha[a,b]$, then we have

$$_aI_b^\alpha f(x) = g(b) - g(a).$$

(2) **(Local fractional integration by parts)** Suppose that $f(x), g(x) \in D_\alpha[a,b]$ and $f^{(\alpha)}(x)$, $g^{(\alpha)}(x) \in C_\alpha[a,b]$, then we have

$$_aI_b^\alpha f(x)g^{(\alpha)}(x) = f(x)g(x)\big|_a^b - {_aI_b^\alpha} f^{(\alpha)}(x)g(x).$$

**LEMMA 2.2.**[10]

$$\frac{d^\alpha x^{k\alpha}}{dx^\alpha} = \frac{\Gamma(1+k\alpha)}{\Gamma(1+(k-1)\alpha)} x^{(k-1)\alpha}; \qquad \frac{1}{\Gamma(1+\alpha)}\int_a^b x^{k\alpha}(dx)^\alpha = \frac{\Gamma(1+k\alpha)}{\Gamma(1+(k+1)\alpha)}(b^{(k+1)\alpha} - a^{(k+1)\alpha}), \ k \in R.$$

**LEMMA 2.3.[10]** (Generalized Hölder's ineauality) Let $f,g \in C_\alpha[a,b]$, $p,q>1$ with $1/p+1/q=1$, then

$$\frac{1}{\Gamma(1+\alpha)}\int_a^b |f(x)g(x)|(dx)^\alpha \leq \left(\frac{1}{\Gamma(1+\alpha)}\int_a^b |f(x)|^p (dx)^\alpha\right)^{1/p}\left(\frac{1}{\Gamma(1+\alpha)}\int_a^b |g(x)|^q (dx)^\alpha\right)^{1/q}.$$

## 3. MAIN RESULTS

**THEOREM 1.** Let $I \subseteq R$ be an interval, $f : I^0 \subseteq R \to R^\alpha$ ($I^0$ is the interior of $I$) such that $f \in D_\alpha(I^\circ)$ and $f^{(\alpha)} \in C_\alpha[a,b]$ for $a,b \in I^0$ with $a<b$. Then the following equality holds:

$$\frac{f(a)+f(b)}{2^\alpha} - \frac{\Gamma(1+\alpha)}{(b-a)^\alpha}\,_aI_b^\alpha f(x) = \frac{(b-a)^\alpha}{2^\alpha}\frac{1}{\Gamma(1+\alpha)}\int_0^1 (1-2t)^\alpha f^{(\alpha)}(ta+(1-t)b)(dt)^\alpha.$$

**Proof:** From the local fractional integration by parts, it suffices to note that

$$\frac{1}{\Gamma(1+\alpha)}\int_0^1 (1-2t)^\alpha f^{(\alpha)}(ta+(1-t)b)(dt)^\alpha$$

$$= (1-2t)^\alpha \frac{f(ta+(1-t)b)}{(a-b)^\alpha}\bigg|_0^1 - \frac{1}{(a-b)^\alpha}\frac{1}{\Gamma(1+\alpha)}\int_0^1 (-2)^\alpha \Gamma(1+\alpha) f(ta+(1-t)b)(dt)^\alpha.$$

$$= \frac{f(a)+f(b)}{(b-a)^\alpha} - \frac{2^\alpha}{(b-a)^{2\alpha}}\frac{\Gamma(1+\alpha)}{\Gamma(1+\alpha)}\int_a^b f(x)(dx)^\alpha$$

$$= \frac{f(a)+f(b)}{(b-a)^\alpha} - \frac{2^\alpha \Gamma(1+\alpha)}{(a-b)^{2\alpha}}\,_aI_b^\alpha f(x).$$

Thus, $\dfrac{(b-a)^\alpha}{2^\alpha}\dfrac{1}{\Gamma(1+\alpha)}\displaystyle\int_0^1 (1-2t)^\alpha f^{(\alpha)}(ta+(1-t)b)(dt)^\alpha = \dfrac{f(a)+f(b)}{2^\alpha} - \dfrac{\Gamma(1+\alpha)}{(b-a)^\alpha}\,_aI_b^\alpha f(x).$

**THEOREM 2.** Suppose that the assumptions of Theorem 1 are satisfied, then

$$f\left(\frac{a+b}{2}\right) - \frac{\Gamma(1+\alpha)}{(b-a)^\alpha}\,_aI_b^\alpha f(x) = \frac{1}{(b-a)^\alpha}\frac{1}{\Gamma(1+\alpha)}\int_a^b S(x)f^{(\alpha)}(x)(dx)^\alpha,$$

where

$$S(x) = \begin{cases}(x-a)^\alpha, & x \in [a,\frac{a+b}{2}], \\ (x-b)^\alpha, & x \in [\frac{a+b}{2},b].\end{cases}$$

**Proof.** From the definition of $S(x)$, we have

$$\frac{1}{(b-a)^\alpha}\frac{1}{\Gamma(1+\alpha)}\int_a^b S(x)f^{(\alpha)}(x)(dx)^\alpha$$

$$= \frac{1}{(b-a)^\alpha} \frac{1}{\Gamma(1+\alpha)} \int_a^{\frac{a+b}{2}} (x-a)^\alpha f^{(\alpha)}(x)(dx)^\alpha + \frac{1}{(b-a)^\alpha} \frac{1}{\Gamma(1+\alpha)} \int_{\frac{a+b}{2}}^b (x-b)^\alpha f^{(\alpha)}(x)(dx)^\alpha.$$ Firstly, let's estimate the first part.

In fact, using the local fractional integration by parts, we have

$$\frac{1}{\Gamma(1+\alpha)} \int_a^{\frac{a+b}{2}} (x-a)^\alpha f^{(\alpha)}(x)(dx)^\alpha = \left(\frac{b-a}{2}\right)^\alpha f\left(\frac{a+b}{2}\right) - \Gamma(1+\alpha) \, _aI^\alpha_{\frac{a+b}{2}} f(x).$$

Secondly, using the local fractional integration by parts, we have

$$\frac{1}{\Gamma(1+\alpha)} \int_{\frac{a+b}{2}}^b (x-b)^\alpha f^{(\alpha)}(x)(dx)^\alpha = \left(\frac{b-a}{2}\right)^\alpha f\left(\frac{a+b}{2}\right) - \Gamma(1+\alpha) \, _{\frac{a+b}{2}}I^\alpha_b f(x).$$

Combining the above formulas, we obtain

$$\frac{1}{(b-a)^\alpha} \frac{1}{\Gamma(1+\alpha)} \int_a^b S(x) f^{(\alpha)}(x)(dx)^\alpha = f\left(\frac{a+b}{2}\right) - \frac{\Gamma(1+\alpha)}{(b-a)^\alpha} \, _aI^\alpha_b f(x).$$

That is

$$f\left(\frac{a+b}{2}\right) - \frac{\Gamma(1+\alpha)}{(b-a)^\alpha} \, _aI^\alpha_b f(x) = \frac{1}{(b-a)^\alpha} \frac{1}{\Gamma(1+\alpha)} \int_a^b S(x) f^{(\alpha)}(x)(dx)^\alpha.$$

**THEOREM 3.** Let $p, q \geq 1$ with $1/p + 1/q = 1$. Let $f : I^0 \subseteq R \to R^\alpha$, such that $f \in D_\alpha(I^\circ)$ and $|f^{(\alpha)}|^q \in C_\alpha[a,b]$ for $a, b \in I^0$ with $a < b$. If $|f^{(\alpha)}|^q$ is generalized convex on $[a,b]$, then

$$\left| \frac{f(a)+f(b)}{2^\alpha} - \frac{\Gamma(1+\alpha)}{(b-a)^\alpha} \, _aI^\alpha_b f(x) \right|$$

$$\leq \frac{(b-a)^\alpha}{2^\alpha} \left[ |f^{(\alpha)}(a)|^q + |f^{(\alpha)}(b)|^q \right]^{1/q} \left( \frac{\Gamma(1+\alpha)}{\Gamma(1+2\alpha)} \right)^{1/p} \left[ \left(\frac{3}{2}\right)^\alpha \frac{\Gamma(1+2\alpha)}{\Gamma(1+3\alpha)} - \left(\frac{1}{2}\right)^\alpha \frac{\Gamma(1+\alpha)}{\Gamma(1+2\alpha)} \right]^{1/q}.$$

**Proof.** From Theorem 1 and the generalized Hölder's inequality (Lemma 3), it follows that

$$\left| \frac{f(a)+f(b)}{2^\alpha} - \frac{\Gamma(1+\alpha)}{(b-a)^\alpha} \, _aI^\alpha_b f(x) \right|$$

$$\leq \frac{(b-a)^\alpha}{2^\alpha} \frac{1}{\Gamma(1+\alpha)} \int_0^1 |(1-2t)^\alpha| \, |f^{(\alpha)}(ta+(1-t)b)| (dt)^\alpha$$

$$\leq \frac{(b-a)^\alpha}{2^\alpha} \left( \frac{1}{\Gamma(1+\alpha)} \int_0^1 |1-2t|^\alpha (dt)^\alpha \right)^{1/p} \left[ \frac{1}{\Gamma(1+\alpha)} \int_0^1 |1-2t|^\alpha |f^{(\alpha)}(ta+(1-t)b))|^q (dt)^\alpha \right]^{1/q}. \quad (3.1)$$

From Lemma 2, it's easy to see that

$$\frac{1}{\Gamma(1+\alpha)} \int_0^1 |1-2t|^\alpha (dt)^\alpha = \frac{\Gamma(1+\alpha)}{\Gamma(1+2\alpha)}. \quad (3.2)$$

Since $|f^{(\alpha)}|^q$ is generalized convex on $[a,b]$, thus

$$\frac{1}{\Gamma(1+\alpha)}\int_0^1|1-2t|^\alpha|f^{(\alpha)}(ta+(1-t)b))|^q\,(dt)^\alpha$$

$$\leq \frac{1}{\Gamma(1+\alpha)}\int_0^1|1-2t|^\alpha[t^\alpha|f^{(\alpha)}(a)|^q+(1-t)^\alpha|f^{(\alpha)}(b)|^q](dt)^\alpha$$

$$=\left[|f^{(\alpha)}(a)|^q+|f^{(\alpha)}(b)|^q\right]\frac{1}{\Gamma(1+\alpha)}\int_0^1|1-2t|^\alpha t^\alpha(dt)^\alpha. \tag{3.3}$$

It is obvious that

$$\frac{1}{\Gamma(1+\alpha)}\int_0^1|1-2t|^\alpha t^\alpha(dt)^\alpha = -\frac{\Gamma(1+\alpha)}{\Gamma(1+2\alpha)}\left(\frac{1}{2}\right)^\alpha + \frac{\Gamma(1+2\alpha)}{\Gamma(1+3\alpha)}\left(\frac{3}{2}\right)^\alpha. \tag{3.4}$$

Combining the estimates of (3.1)-(3.4), we have

$$\left|\frac{f(a)+f(b)}{2^\alpha}-\frac{\Gamma(1+\alpha)}{(b-a)^\alpha}\,_aI_b^\alpha f(x)\right|$$

$$\leq \frac{(b-a)^\alpha}{2^\alpha}\left[|f^{(\alpha)}(a)|^q+|f^{(\alpha)}(b)|^q\right]^{1/q}\left(\frac{\Gamma(1+\alpha)}{\Gamma(1+2\alpha)}\right)^{1/p}\left[\frac{\Gamma(1+2\alpha)}{\Gamma(1+3\alpha)}\left(\frac{3}{2}\right)^\alpha - \frac{\Gamma(1+\alpha)}{\Gamma(1+2\alpha)}\left(\frac{1}{2}\right)^\alpha\right]^{1/q}.$$

**THEOREM 4.** Suppose that the assumptions of Theorem 3 are satisfied, then

$$\left|f\left(\frac{a+b}{2}\right)-\frac{\Gamma(1+\alpha)}{(b-a)^\alpha}\,_aI_b^\alpha f(x)\right|\leq \frac{(b-a)^\alpha}{4^\alpha}\left[\frac{2^\alpha\Gamma(1+\alpha)}{\Gamma(1+2\alpha)}\right]^{1/p}\left(\frac{\Gamma(1+\alpha)[|f^{(\alpha)}(a)|^q+|f^{(\alpha)}(b)|^q]}{\Gamma(1+2\alpha)}\right)^{1/q}.$$

**Proof.** Similar to the estimate of Theorem 3, from Theorem 2 and Lemma 3, it follows that

$$\left|f\left(\frac{a+b}{2}\right)-\frac{\Gamma(1+\alpha)}{(b-a)^\alpha}\,_aI_b^\alpha f(x)\right|$$

$$\leq \frac{1}{(b-a)^\alpha}\frac{1}{\Gamma(1+\alpha)}\int_a^b|S(x)||f^{(\alpha)}(x)|(dx)^\alpha$$

$$\leq \frac{1}{(b-a)^\alpha}\left[\frac{1}{\Gamma(1+\alpha)}\int_a^b|S(x)|(dx)^\alpha\right]^{1/p}\left[\frac{1}{\Gamma(1+\alpha)}\int_a^b|S(x)||f^{(\alpha)}(x)|^q(dx)^\alpha\right]^{1/q}.$$

Let us estimate the first integral.

Using the expression of function $S(x)$, we can calculate that

$$\frac{1}{\Gamma(1+\alpha)}\int_a^b|S(x)|(dx)^\alpha = \frac{\Gamma(1+\alpha)}{\Gamma(1+2\alpha)}2^\alpha\left(\frac{b-a}{2}\right)^{2\alpha}. \tag{3.5}$$

Secondly, for the another part, we have

$$\frac{1}{\Gamma(1+\alpha)}\int_a^b|S(x)||f^{(\alpha)}(x)|^q(dx)^\alpha$$

$$= \frac{1}{\Gamma(1+\alpha)} \int_a^{\frac{a+b}{2}} (x-a)^\alpha \mid f^{(\alpha)}(x) \mid^q (dx)^\alpha + \frac{1}{\Gamma(1+\alpha)} \int_{\frac{a+b}{2}}^b (b-x)^\alpha \mid f^{(\alpha)}(x) \mid^q (dx)^\alpha.$$

Since $\mid f^{(\alpha)} \mid^q$ is generalized convex on $[a,b]$, thus we have

$$\frac{1}{\Gamma(1+\alpha)} \int_a^{\frac{a+b}{2}} (x-a)^\alpha \mid f^{(\alpha)}(x) \mid^q (dx)^\alpha$$

$$\leq \frac{1}{\Gamma(1+\alpha)} \int_0^1 \left(\frac{b-a}{2}\right)^{2\alpha} t^\alpha \left| f^{(\alpha)}\left(\frac{a+b}{2}t + (1-t)a\right) \right|^q (dt)^\alpha$$

$$\leq \frac{1}{\Gamma(1+\alpha)} \int_0^1 \left(\frac{b-a}{2}\right)^{2\alpha} t^\alpha \left[ t^\alpha \left| f^{(\alpha)}\left(\frac{a+b}{2}\right) \right|^q + (1-t)^\alpha \mid f^{(\alpha)}(a) \mid^q \right] (dt)^\alpha$$

$$\leq \left(\frac{b-a}{2}\right)^{2\alpha} \left\{ \frac{\mid f^{(\alpha)}(a) \mid^q + \mid f^{(\alpha)}(b) \mid^q}{2^\alpha} \frac{\Gamma(1+2\alpha)}{\Gamma(1+3\alpha)} + \left| f^{(\alpha)}(a) \right|^q \left[ \frac{\Gamma(1+\alpha)}{\Gamma(1+2\alpha)} - \frac{\Gamma(1+2\alpha)}{\Gamma(1+3\alpha)} \right] \right\}. \quad (3.6)$$

Similarly, we have

$$\frac{1}{\Gamma(1+\alpha)} \int_{\frac{a+b}{2}}^b (b-x)^\alpha \mid f^{(\alpha)}(x) \mid^q (dx)^\alpha$$

$$\leq \left(\frac{b-a}{2}\right)^{2\alpha} \left\{ \frac{\mid f^{(\alpha)}(a) \mid^q + \mid f^{(\alpha)}(b) \mid^q}{2^\alpha} \frac{\Gamma(1+2\alpha)}{\Gamma(1+3\alpha)} + \left| f^{(\alpha)}(b) \right|^q \left[ \frac{\Gamma(1+\alpha)}{\Gamma(1+2\alpha)} - \frac{\Gamma(1+2\alpha)}{\Gamma(1+3\alpha)} \right] \right\}. \quad (3.7)$$

So, from (3.6) and (3.7), we get

$$\frac{1}{\Gamma(1+\alpha)} \int_a^b \mid S(x) \mid \mid f^{(\alpha)}(x) \mid^q (dx)^\alpha = \left(\frac{b-a}{2}\right)^{2\alpha} \left( \left| f^{(\alpha)}(a) \right|^q + \left| f^{(\alpha)}(b) \right|^q \right) \frac{\Gamma(1+\alpha)}{\Gamma(1+2\alpha)}. \quad (3.8)$$

Thus, combining (3.5) and (3.8), we obtain

$$\left| f\left(\frac{a+b}{2}\right) - \frac{\Gamma(1+\alpha)}{(b-a)^\alpha} {}_aI_b^\alpha f(x) \right| = \frac{(b-a)^\alpha}{4^\alpha} \left[ \frac{2^\alpha \Gamma(1+\alpha)}{\Gamma(1+2\alpha)} \right]^{1/p} \left( \frac{\Gamma(1+\alpha)[\mid f^{(\alpha)}(a) \mid^q + \mid f^{(\alpha)}(b) \mid^q]}{\Gamma(1+2\alpha)} \right)^{1/q}.$$

## 4. APPLICATIONS TO SOME SPECIAL MEANS

As in [5], in the paper we consider some generalized means. We take

$$A(a,b) = \frac{a^\alpha + b^\alpha}{2^\alpha}; \qquad L_n(a,b) = \left[ \frac{\Gamma(1+n\alpha)}{\Gamma(1+(n+1)\alpha)} (b^{(n+1)\alpha} - a^{(n+1)\alpha}) \right]^{1/n}, \quad n \in Z \setminus \{-1, 0\}, \quad a,b \in R, \quad a \neq b.$$

**Proposition 1.** Let $a, b \in R$, $a < b$, $0 \notin [a,b]$ and $n \in Z$, $\mid n \mid \geq 2$. Then for all $p, q \geq 1$ with $1/p + 1/q = 1$,

$$\mid A(a^n, b^n) - \frac{\Gamma(1+\alpha)}{(b-a)^\alpha} (L_n(a,b))^n \mid$$

$$\leq \frac{(b-a)^\alpha}{2^\alpha} \left[ \left( \frac{\Gamma(1+n\alpha)}{\Gamma(1+(n-1)\alpha)} \right)^q (a^{(n-1)q\alpha} + b^{(n-1)q\alpha}) \right]^{1/q} \left( \frac{\Gamma(1+\alpha)}{\Gamma(1+2\alpha)} \right)^{1/p} \left[ \frac{\Gamma(1+2\alpha)}{\Gamma(1+3\alpha)} \left(\frac{3}{2}\right)^\alpha - \frac{\Gamma(1+\alpha)}{\Gamma(1+2\alpha)} \left(\frac{1}{2}\right)^\alpha \right]^{1/q};$$

$$\left| (A(a,b))^n - \frac{\Gamma(1+\alpha)}{(b-a)^\alpha}(L_n(a,b))^n \right| \leq \frac{(b-a)^\alpha}{4^\alpha}\left[ \frac{2^\alpha \Gamma(1+\alpha)}{\Gamma(1+2\alpha)} \right]^{1/p} \left[ \frac{\Gamma(1+\alpha)}{\Gamma(1+2\alpha)} \left( \frac{\Gamma(1+n\alpha)}{\Gamma(1+(n-1)\alpha)} \right)^q (a^{(n-1)q\alpha} + b^{(n-1)q\alpha}) \right]^{1/q}.$$

**Proof.** The proof is from Theorem 3 and Theorem 4 with $f(x) = x^{n\alpha}$, $x \in R$, $n \in Z$, $n \geq 2$.

$$|A(a^n, b^n) - \frac{\Gamma(1+\alpha)}{(b-a)^\alpha}(L_n(a,b))^n| = \left| \frac{f(a)+f(b)}{2^\alpha} - \frac{\Gamma(1+\alpha)}{(b-a)^\alpha} {}_aI_b^\alpha f(x) \right|$$

$$\leq \frac{(b-a)^\alpha}{2^\alpha}\left[ \left( \frac{\Gamma(1+n\alpha)}{\Gamma(1+(n-1)\alpha)} \right)^q (a^{(n-1)q\alpha} + b^{(n-1)q\alpha}) \right]^{1/q} \left( \frac{\Gamma(1+\alpha)}{\Gamma(1+2\alpha)} \right)^{1/p} \left[ \frac{\Gamma(1+2\alpha)}{\Gamma(1+3\alpha)}\left(\frac{3}{2}\right)^\alpha - \frac{\Gamma(1+\alpha)}{\Gamma(1+2\alpha)}\left(\frac{1}{2}\right)^\alpha \right]^{1/q}.$$

And,

$$|(A(a,b))^n - \frac{\Gamma(1+\alpha)}{(b-a)^\alpha}(L_n(a,b))^n| = \left| f\left(\frac{a+b}{2}\right) - \frac{\Gamma(1+\alpha)}{(b-a)^\alpha} {}_aI_b^\alpha(x) \right|$$

$$\leq \frac{(b-a)^\alpha}{4^\alpha}\left[ \frac{2^\alpha \Gamma(1+\alpha)}{\Gamma(1+2\alpha)} \right]^{1/p} \left[ \frac{\Gamma(1+\alpha)}{\Gamma(1+2\alpha)} \left( \frac{\Gamma(1+n\alpha)}{\Gamma(1+(n-1)\alpha)} \right)^q (a^{(n-1)q\alpha} + b^{(n-1)q\alpha}) \right]^{1/q}.$$

**Remark:** when $\alpha = 1$, this result reduce to Proposition 1 in [5].

## ACKNOWLEDGMENTS

This work has been supported by the National Natural Science Foundation of China (No. 11161042, 11471050).